\documentclass[conference]{IEEEtran}
\pdfoutput=1
\usepackage{keyval}
\AtBeginDocument{\usepackage{graphicx}}

\usepackage{color}
\usepackage{epsfig}
\usepackage{scalefnt}
\newcommand {\bel}[1]{\begin{align*}}
\newcommand {\eel}[1]{\end{align*}}
\newcommand {\bea}{\begin{eqnarray}}
\newcommand {\eea}{\end{eqnarray}}

\newcommand{\ignore}[1]{\relax}

\renewcommand{\epsilon}{\varepsilon}

\newcommand{\comments}[1]{}

\title{Locating PHEV Exchange Stations in V2G}
\author{\IEEEauthorblockN{Feng Pan, Russell Bent, 
Alan Berscheid,  David Izraelevitz }
\IEEEauthorblockA{ 
Decision Applications Division\\
Los Alamos National Laboratory
}
}

\begin{document}
\maketitle 
\begin{abstract}
Plug-in hybrid electric vehicles (PHEVs) are an environmentally friendly technology that is expected to rapidly penetrate the transportation system. Renewable energy sources such as wind and solar have received considerable attention as clean power options for future generation expansion. However, these sources are intermittent and increase the uncertainty in the ability to generate power.  The deployment of PHEVs in a vehicle-to-grid (V2G) system provide a potential mechanism for reducing the variability of renewable energy sources. For example, PHEV supporting infrastructures like battery exchange stations that provide battery service to PHEV customers could be used as storage devices to stabilize the grid when renewable energy production is fluctuating. In this paper, we study how to best site these stations in terms of how they can support both the transportation system and the power grid. To model this problem  we develop a two-stage stochastic program to optimally locate the stations prior to the realization of battery demands, loads, and generation capacity of renewable power sources. We develop two test cases to study the benefits and the performance of these systems.

\end{abstract}
\vspace{-.2cm}
\section*{Nomenclature}
\vspace{-.2cm}

\addcontentsline{toc}{section}{Nomenclature}
\noindent{\bf Index/Set} 
\begin{IEEEdescription}[\IEEEusemathlabelsep\IEEEsetlabelwidth{$V_1,V_2,V_3$}]
\item[$G(N,E)$] power grid with bus set $N$ and link set $E$
\item [$J$] set of traffic routes
\item [$I$] set of exchange stations
\item [$u = m(i)$] the bus, $u$, of station $i$
\item [$\Omega$] set of scenarios $\{\omega_1,\omega_2,\ldots\}$
\end{IEEEdescription}
{\bf Deterministic Data}
\begin{IEEEdescription}[\IEEEusemathlabelsep\IEEEsetlabelwidth{$V_1,V_2,V_3$}]
%{\setlength{\leftmargin}{3cm}}
\item [$f_i$] fixed cost to open an exchange station at $i$
\item [$r_i$] storage cost per battery at station $i$
\item [$U_i$] maximum number of batteries at station $i$
\item [$L_i$] minimum number of batteries at station $i$
\item [$c_{ij}$] cost to use station $i$ from traffic route $j$
\item [$h_j$] cost of unmet demand from traffic route $j$ 
\item [$a$] power per battery (MW)
\item [$C_{uv}$] line capacity of line $(u,v)$
\item [$g_u$] load shedding cost at bus $u$  
\item [$b_{uv}$] reactance of line $(u,v)$
\end{IEEEdescription}
{\bf Stochastic Data (per scenario $\omega$}
\begin{IEEEdescription}[\IEEEusemathlabelsep\IEEEsetlabelwidth{$V_1,V_2,V_3$}]
%{\setlength{\leftmargin}{3cm}}
\item [$d^\omega_j$] demand for batteries from traffic route $j$
\item [$l^\omega_{u}$] load at bus $u$
\item [$G^\omega_{u}$] generation capacity at $u$
\item [$o^\omega_{u}$] generation cost per MW at $u$
\end{IEEEdescription}
{\bf Variables}
\begin{IEEEdescription}[\IEEEusemathlabelsep\IEEEsetlabelwidth{$V_1,V_2,V_3$}]
%{\setlength{\leftmargin}{3cm}}
\item [$x_i$] 0-1 variable siting an exchange station at $i$
\item [$w_i$] number of batteries stored at station $i$
\item [$t^\omega_i$] number of batteries for PHEV demand
\item [$s^\omega_j$] number of batteries discharged to the grid
\item [$q^\omega_i$] number of unserved batteries on route $j$
\item [$y^\omega_{ij}$] number of batteries served at station $i$ from route $j$
%%\item [$r^\omega_{j}$] number of unserved batteries from traffic pattern $j$
\item [$\alpha^\omega_{uv}$] power flow in line $(u,v)$
\item [$\theta^\omega_{u}$] phase angle at bus $u$
\item [$\delta^\omega_{u}$] slack variable for unmet demand at bus $u$
\item [$\beta^\omega_{u}$] power generated at bus $u$
\end{IEEEdescription}

\vspace{-.2cm}

\section{Introduction}
\begin{figure}
\centering
\psfig{file=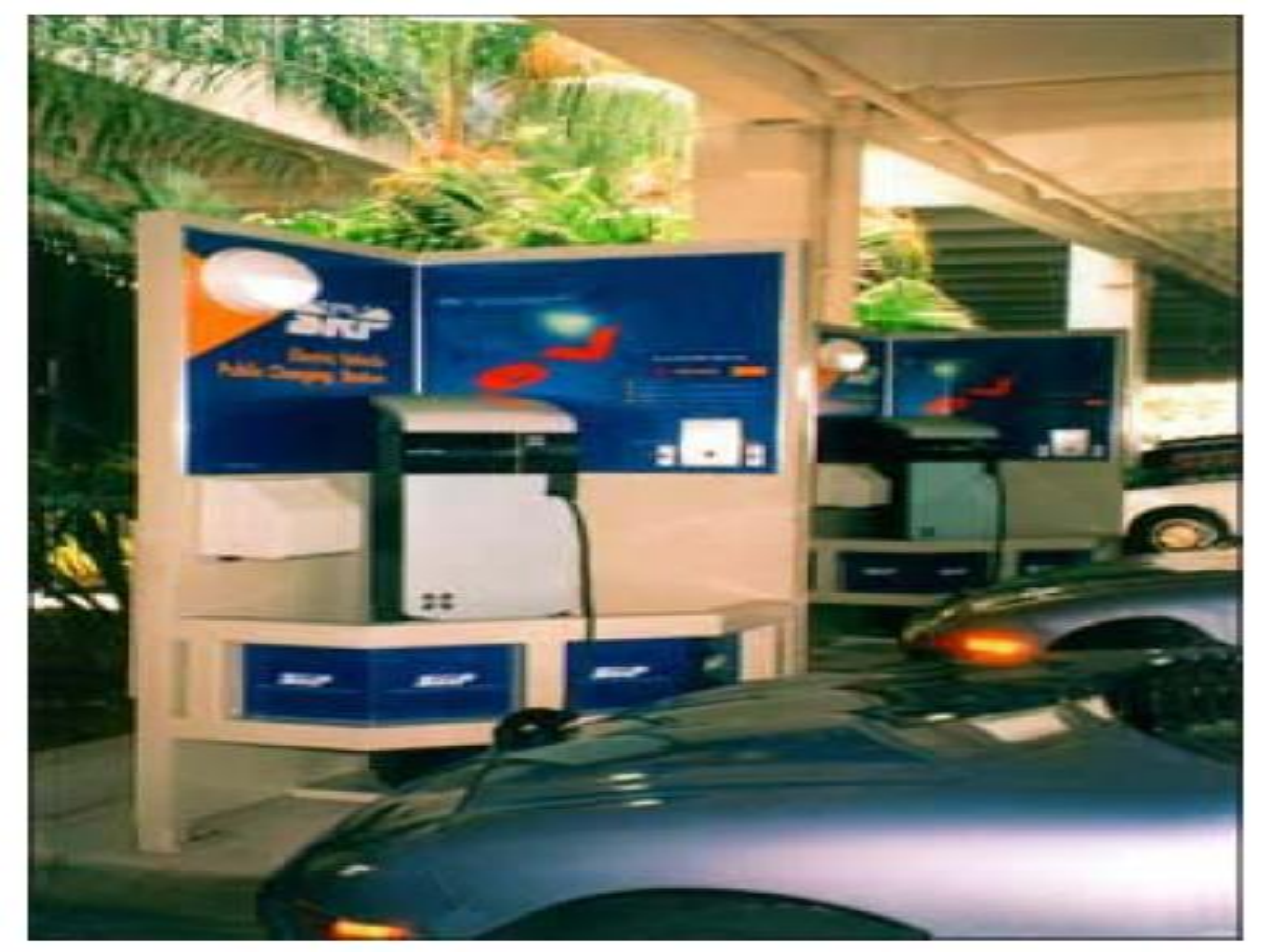,scale=.16,page=1} \hspace{.2cm}
\psfig{file=PHEV-Figures.pdf,scale=.16,page=2} \hspace{.2cm}
\caption{\label{fig:chargingstations} \small Commercial charging stations }
\vspace{-0.2cm}
\end{figure}

\begin{figure}
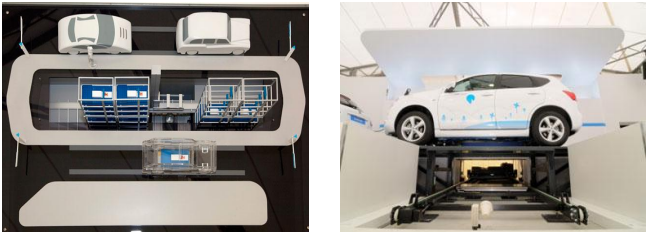

\centering
\psfig{file=PHEV-Figures.pdf,scale=.16,page=3} \hspace{.2cm}
\psfig{file=PHEV-Figures.pdf,scale=.16,page=4} \hspace{.2cm}
\caption{\label{fig:exchangestations} \small Commercial exchange stations }
\vspace{-0.6cm}
\end{figure}

%%--importance of phev 
Global warming and dependence on fossil fuels pose great challenges to the nation's energy infrastructures and future energy consumption. Smart grid technologies and renewable generation are often touted as central pieces to address these challenges \cite{09-DOE-smartgrid}. Indeed, there have been numerous recent studies analyzing the cost-benefits of large penetration of renewable generation into the power grid \cite{07-Boy, 07-Flo, 09-Too} and the implementation of vehicle-to-grid (V2G) \cite{01-Kem,09-Nem} systems for PHEV to reduce fossil fuel consumption \cite{07-web-history-hv,10-DOE-Veh,07-Gon-phev}. 
%Modern hybrid electric vehicles were introduced into the automotive market in 1997 \cite{07-web-history-hv}.  As of 2009, there were 1.6 million PHEVs vehicles in the U.S. \cite{10-DOE-Veh}
%and this number has increased at the rate of 80\% per year \cite{07-Gon-phev}. 
Recent studies \cite{08-Mor-phev, 06-Sim-phev,08-DOT} on feasibility, cost-and-benefit, and impact of PHEVs indicate that it is necessary to build supporting PHEV infrastructures (charging stations (see Fig.~\ref{fig:chargingstations} and \cite{08-Mor-phev}) and exchange stations (see Fig.~\ref{fig:exchangestations} and \cite{exchangestation})) if the adoption of hybrid electric vehicles continues to increase at current rates. There are challenges in implementing both systems: variability in renewable generation and increased demand for power. In this paper, we focus on how renewable technologies and V2G exchange stations may be coupled to address these challenges.

The crucial connection between the two systems arises from the observation that exchange stations are essentially large battery banks where PHEV drivers can automatically switch their batteries with fully charged batteries.  The exchange station may choose to charge batteries during low power usage times (off-peak) or when renewable power plants are producing large amounts of power, thereby reducing the demands a V2G system places on the grid.  In addition, the station may discharge batteries on to the grid during periods of low renewable generation thereby reducing the variability of renewable generation. Thus, there is considerable potential that a coupling of the two systems will reduce the challenges that arise when the systems are implemented independently. This paper seeks to study how to site the stations so that they benefit both the grid and the V2G system.
This problem shares many similarities to the facility location problem \cite{95-Dre,04-Sny,04-Eis, 06-Wu,02-Dre}. In the facility location problem,  facilities are opened to serve a set of customer demands, and the objective is to minimize the setup cost of opening facilities and transportation cost for customers to use the facility. 
%There are many variations of the facility location problem, e.g., facility with capacity \cite{04-Eis, 06-Wu} and the stochastic customer demands \cite{08-Sch,04-Sny}. 
%For more comprehensive review of facility location problems, see \cite{95-Dre, 02-Dre, 04-Sny}. 
The crucial difference between the problem stated in this paper and traditional facility location is that this problem problem contains two networks, a transportation network and a power grid. 
%Considering the power grid,  the exchange stations can be viewed as battery banks which provide the addition generation capacity. Thus considering the power grid alone,  we have  the generation expansion problem which is important for planning the grid and has been studied extensively \cite{07-Lop,97-Zhu, 07-Kay}. Generation expansion problem minimize the overall cost to allocation power generation capacity in a grid to satisfy the demands.  Some variants of the problem also incorporate stochastic  demand  \cite{10-Gar} and non-regulated markets \cite{06-Ehr}.  

The contributions of this paper are two-fold.  First, we develop the first mathematical model to site exchange stations that accounts for impacts to transportation systems and electric power systems. Second, we investigate variations of this model to determine the necessity of accounting for both systems, the conditions in which proper siting can be beneficial to both systems, and potential benefit s trade-offs  to both systems.
The rest of the paper is outlined as follows. 
%Charging stations and exchange stations are discussed in  Section~\ref{sec:stations}.   
A two-stage stochastic program is developed for designing the location of exchange stations in a V2G in Section~\ref{sec:model}. In Section~\ref{sec:experiment}, we investigate the effect of exchange stations on power grids with renewable generation and we conclude with Section~\ref{sec:conclusion}. 

%\section{PHEV exchange stations and charging stations} \label{sec:stations}

\vspace{-.2cm}

\section{A Two-stage Stochastic Program for Siting Exchange Stations}\label{sec:model}
In a V2G system, exchange stations can serve PHEVs and provide storage services to reduce the variability introduced by renewable generation sources. Thus, it is strategically important to site  exchange stations in locations that are convenient for PHEV drivers and that can be connected to the grid without requiring additional capacity expansion. %, e.g., putting additional power lines to accommodate the power exchange between the grid and  stations.ent for PHEV drivers and that can be connected to the grid without requiring additional capacity expansion in a grid.%, e.g., putting additional power lines to accommodate the power exchange between the grid and  stations.
In this section, we introduce a two-stage stochastic program that integrates a battery demand model and a power grid model to optimally site exchange stations in a V2G system.

In our stochastic programming model, a time line is defined over the decision variables and scenarios as follows. The (first-stage) decisions on where to site exchanges stations and how many batteries to store at each station are made prior to knowing the future demand for the batteries (scenarios). The scenarios are drawn from discrete random variables that characterize PHEV battery demands on each traffic route, the load at buses, and production capacity of renewable generators.  After a scenario is realized, recourse (second-stage) decisions are made to distribute batteries to satisfy PHEV demand and discharge power back to the grid. The cost of the recourse decisions is dependent on the first-stage decisions and the random scenario.  The overall objective function combines the first-stage cost and the expected cost of recourse actions over all scenarios.  

We now formalize the model of the first stage variables. In the first stage, the location and size of exchange stations are decided. The long term decisions are made without knowing the exact PHEV demands, loads, and generation capacities. The objective function is to minimize the overall cost
\begin{equation}
\min_{x,w}   \sum_{i\in I}{(f_ix_i+r_i w_i)} + E_\Omega[h(x,w,\xi)] \label{eqn:obj1}
\end{equation}
In (\ref{eqn:obj1}), the summation includes $f_ix_i$  and $r_iw_i$ which represents the fixed cost to open exchange station $i$ and the cost to store batteries at $i$. A first-stage  constraint is represented by the following equation:
\begin{equation}
l_i x_i \le w_i \le u_i x_i, \hspace{5pt} \forall i \in I \label{eqn:cap_batteries}.
\end{equation}
This constraint states that each exchange station can only store batteries when it is open and must satisfy a lower and upper bound on the storage capacity. 
The expectation term in the objective represents the second-stage recourse cost $h(x,w,\xi) $  of satisfying PHEV demands and meeting demands for power over a set $\Omega$ of scenarios. 

A scenario specifies a realization of the random variables.  For each scenario, $\omega$, the second stage objective $h(x,w,\xi)$ is defined as
\begin{equation}
\min \sum_{i\in I, j\in J}{c_{ij}y^\omega_{ij}} +\sum_{j\in J}{h^\omega_j q^\omega_j}+\sum_{u\in N}{o^\omega_{u} \beta^\omega_{u}} + \sum_{u\in N}{g_u \delta^\omega_u} \label{eqn:obj2}.
\end{equation}
%This objective include the both costs for evaluating the PHEV battery service and the power grid.
The term $c_{ij}y^\omega_{ij}$ represents the cost of PHEV drivers on route $j$ to swap batteries at station $i$. The penalty cost on unmet PHEV battery demand from route $j$ is stated with $h^\omega_jq^\omega_j$. The cost $c_{ij}$ from route $j$ to station $i$ is typically a function of the distance a driver must divert from route $j$ to reach station $i$. For the power grid, the costs include $o^\omega_{u} \beta^\omega_{u}$ for generation and $g_u \delta^\omega_u$ for shedding loads, which are used to measure the performance of the power grid and include the availability of batteries to meet demand for power. 

At each station, available batteries $w_i$ can be used to satisfy PHEV demand ($t^\omega_i$) or supply power to the grid ($s^\omega_i$). The following constraint is used
to ensure that no station uses more batteries than are stored 
\begin{equation}
s^\omega_i + t^\omega_i \le w_i,  \hspace{5pt} i \in I \label{eqn:control_station}.
\end{equation}

For the transportation system in the PHEV battery model, we model traffic as routes. The traffic density at each route is derived from activities which are based on demographic surveys and activity surveys collected from real households in the study areas.  This approach has been used to simulate metropolitan area traffic and a detailed description of the traffic model is found in \cite{FastTrans}.   For each route $j$, there is a demand $d^\omega_j$ for PHEV batteries. In (\ref{eqn:bal_phev}), the demand can be satisfied by a battery at a station or be unsatisfied with some penalty, 
\begin{equation}
 \sum_{i\in I}{y^\omega_{ij}}  + q^\omega_j =  d^\omega_j, \hspace{5pt}.  j \in J \label{eqn:bal_phev}.
\end{equation}
Constraint (\ref{eqn:cap_station}) limits the number batteries used to satisfy PHEV demand to the number of batteries allocated to $t^\omega_i$ at each station.  
\begin{equation}
\sum_{j\in J}{y^\omega_{ij}} \le t^\omega_i, \hspace{5pt}  i \in I \label{eqn:cap_station}
\end{equation}

We use the so called linearized DC power flow equations for modeling flow in the power grid. Constraints (\ref{eqn:bal_bus}) and (\ref{eqn:cap_line}) are the standard constraints for capacitated network problem: (\ref{eqn:bal_bus}) is the network flow balance constraint at each node and (\ref{eqn:cap_line}) is the capacity constraint for each line. In (\ref{eqn:bal_bus}), the  extra term $\sum_{i:m(i) = u}{a s^\omega_i}$ accounts for the power supplied by  the batteries at station $i$ which is connected to bus $u$. 
%In the setting of power grid, 
%power flow can move in either direction in a line, and the line capacity  is proportional to the  line thermal limit. 
In this model, constraint (\ref{eqn:flow_phase}) captures the relationship between power flow on a line and the phase angle difference at either end of the line. 
%This constraint is often used to linearly approximate the non-linear power equations.
At each generator, the generation capacity is enforced by the constraint (\ref{eqn:cap_generation}). For renewable generators, the capacities may vary between scenarios according to a distribution governing the potential output of the generators. For non-renewable generators, the capacity is fixed throughout the scenarios. 

\begin{equation}
\sum_{v: (u,v) \in E}{\alpha^\omega_{uv}} = -l^\omega_u + \delta^\omega_{u} + \beta^\omega_u + \sum_{i:m(i) = u}{a s^\omega_i}, \hspace{5pt}  \forall u\in N \label{eqn:bal_bus}
\end{equation}
\begin{equation}
\alpha^\omega_{uv} = (\theta^\omega_u - \theta^\omega_v)/b_{uv}, \hspace{5pt} \forall (u,v) \in E \label{eqn:flow_phase}
\end{equation}
\begin{equation}
 -C_{uv} \le \alpha^\omega_{uv} \le  C_{uv}, \hspace{5pt} \forall (u,v) \in E \label{eqn:cap_line}
\end{equation}
\begin{equation}
\beta^\omega_{u} \le G^\omega_u, \hspace{5pt} \forall u \in N \label{eqn:cap_generation}
\end{equation}

The equations (\ref{eqn:obj1})--(\ref{eqn:cap_generation}) define a two-stage stochastic program for locating PHEV exchange stations in the V2G system, and we denote this model as {\bf V2G-2STAGES}.

%\subsection{Algorithm}
%--------------- need to add more or remove it ---------------

%This two-stage model consists two optimization models. If we remove all the constraints and objectives related to power grid model, 
%the second stage is an assignment problem and the entire optimization becomes a stochastic facility location problem. On the other 
%hand, if we remove all the constraints and objectives related to PHEV, the second stage becomes a DC optimal flow model and the 
%two-stage model is a stochastic generation expansion. This  two-stage model is NP-Complete, and  we can separate the second stage 
%problem by relaxing (\ref{eqn:control_station}).

\vspace{-.2cm}

\section{Computational Experiments}\label{sec:experiment}

To evaluate the impact of exchange stations in a V2G system and the power grid, we test {\bf V2G-2STAGES} on a variation of the IEEE RTS-79 \cite{79-RTS} benchmark and a problem derived from power grid and transportation datasets for Miami, Florida. Before describing the details of the two case studies, we first discuss the general settings and intended findings for both studies. 
Load shedding and unmet battery demands are two measures of the performance of the V2G system. Although the {\bf V2G-2STAGES} model optimally allocates exchange stations to meet battery demands and stabilize a grid with renewable energy resources, it is important to realize that the primary purpose of exchange stations is to serve PHEVs.  In our case studies, we will investigate the trade-off between these two objectives. Through these studies, we aim to answer the following questions: {\em 1) What is the impact of exchange stations in a V2G system that is connected to a power grid with variable renewable generation, 2) What are the trade-offs between the two objectives, and 3) Is it important to consider both networks for strategic planning of exchange station siting? } 

To model renewable generators, we vary renewable penetration from 0 to 1 in increments of 0.1. For example at renewable penetration level 0.3, each generator may become a renewable source with probability 0.3. The PHEV battery demand is derived from total population, vehicle-to-population ratio (0.78), phev-to-vehicle ratio (0.1), and 10\% of PHEVs requring battery exchanges in any given scenario.  The demand is allocated to each traffic route weighted by the route's average utilization which is set to be same among all routes in the studies.
One hundred scenarios are generated for the second stage. The random variables include battery demand from PHEVs, demand for power, and generation capacities of renewables. All random variables are assumed to be independent (though in the future we plan to explore dependencies in renewable generation capacities). For a given renewable generator, its generation capacity can be 0, 0.5, or 1 of its maximum generation capacity. Load at a bus is a uniform random variable between 0.5 and 1.0 of its peak load. PHEVs demand on a traffic route is a uniform random variable between 0.5 and 1.5 of the average battery demand.

 For both case studies, we implemented {\bf V2G-2STAGES} in C++. The resulting mixed integer program is solved by the branch-and-bound algorithm in CPLEX 11.0 with an optimality gap of 0.01. The longest computing time was under a minute on a standard desktop personal PC.

\subsection{Case Study: IEEE RTS-79}
\begin{figure}
\centering
\psfig{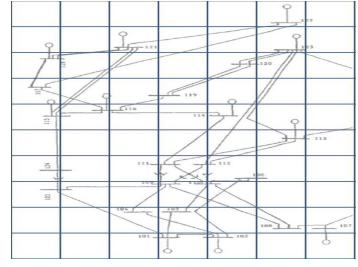} \hspace{.25cm}
\caption{\label{fig:2grids} \small Transportation grid and power grid}
\vspace{-0.6cm}
\end{figure}

Our first case study considers a synthetic city based on the IEEE RTS-79 benchmark designed to mimic the structure of Los Alamos, New Mexico. In this data set, there are  25 buses  and 38 power lines. Loads are distributed at 24 buses and there are 11  buses which have generation capacity, with up to 6 generators located at a single bus. The maximum load and generation capacity are provided for every bus in the data set. The total generation capacity is 2999 MW, and the total demand is 2880 MW. 
We used an 8 by 11 lattice as the transportation network and created the synthetic city (Fig.~\ref{fig:2grids}) by connecting the lattice and power grid. The connection between the transportation network and the power grid was done by assigning the nearest bus to each of the 88 nodes in the transportation network.  
The population for the synthetic city is set to 344,850 which was derived from the population-load ratio of Los Alamos, New Mexico.  
Twenty-eight nodes were randomly selected from the lattice as possible locations for exchange stations. Ten traffic routes were created by randomly selecting ten origin-destination pairs and finding a shortest path for each pair. This synthetic city is a model  of a city with some power consuming industries and a relatively low population. 

%Twenty-eight nodes were randomly selected from the lattice as possible locations for exchange stations. Ten traffic routes were created by randomly selecting ten origin-destination pairs and finding a shortest path for each pair. This synthetic city is a model  of a city with some power consuming industries and a relatively low population. 

\begin{figure}
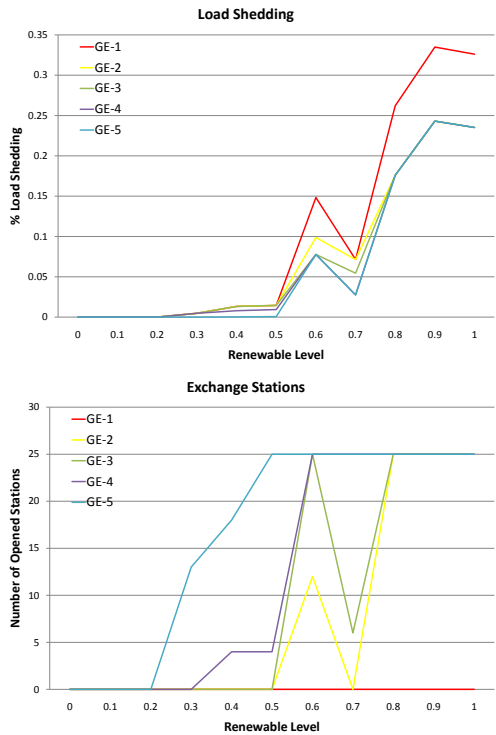

\centering
\psfig{file=PHEV-Figures.pdf,scale=.26,page=11} 
\psfig{file=PHEV-Figures.pdf,scale=.26,page=12}
\caption{\label{fig:rts-ge} \small RTS-79 load shedding and unmet battery demands in generation expansion model with renewables from 0 to 100 percent}
\vspace{-0.6cm}
\end{figure}
In order to understand how V2G batteries may maximally improve the performance of a power grid with variable renewable generation, we
removed PHEV demands and reduced {\bf V2G-2STAGES} to a generation expansion model.  In the base case (GE-1 in Fig.~\ref{fig:rts-ge}), no additional generation resources (exchange stations) were added and the grid was evaluated with different renewable penetration levels. With renewable levels lower than 30\%, there is no loading shedding due to renewable generation variability. For higher renewable penetration level, load shedding can reach as high as 33\%. In cases GE-2 to GE-5, we varied the load shedding penalty to stimulate the opening of exchange stations as generation capacity reserve. We ran all cases with 11 different renewable penetration levels from 0 to 100 percent.  When load shedding penalities are increased, there is some reduction in load shedding. Load shedding is reduced to below 10\% with renewable penetration levels less than 80\%. The reduction is mostly caused by opening 25 out of 28 stations. The opening of additional stations does not reduce load shedding due to other limits, e.g., line capacity, of the grid. 
\begin{figure}
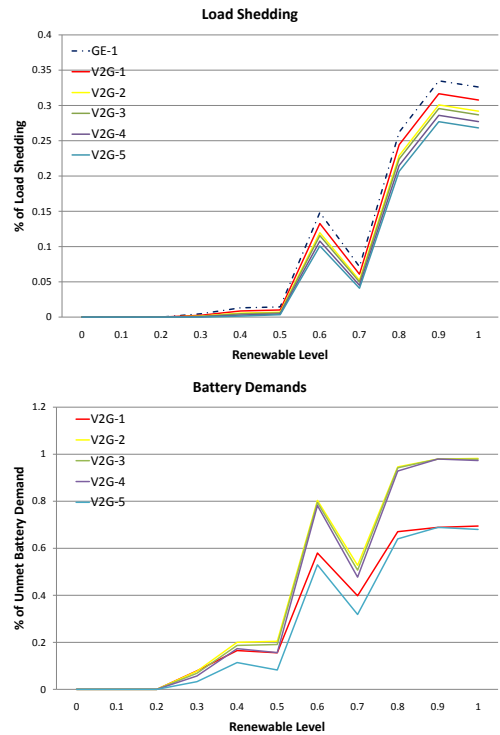

\centering
\psfig{file=PHEV-Figures.pdf,scale=.26,page=13} 
\psfig{file=PHEV-Figures.pdf,scale=.26,page=14}
\caption{\label{fig:rts-v2g1} \small RTS-79: load shedding and unmet battery demands in {\bf V2G-2STAGES} with renewables from 0 to 100 percent}
\vspace{-0.6cm}
\end{figure}

\begin{figure}
\centering
\psfig{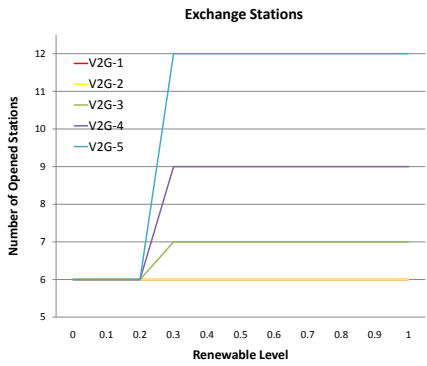} 
\caption{\label{fig:rts-v2g2} \small RTS-79: allocation of exchange stations in {\bf V2G-2STAGES} with renewables from 0 to 100 percent}
\vspace{-0.6cm}
\end{figure}

Fig.~\ref{fig:rts-v2g1} shows the simulations of {\bf V2G-2STAGES} for the entire V2G system (including both PHEV battery demands and power grid). First, we used {\bf V2G-2STAGES} to determine that the optimal solution  to satisfy all battery demands  is 6 stations. In simulation V2G-1, we fixed the station locations to this optimal solution and added the power grid.   There was some reduction in load shedding, but there is a trade-off in unsatisfied PHEV demand. In simulation V2G-2 to V2G-5, we relaxed the location of the stations and allowed the total number of opened stations to be 6, 8, 10, and 12 stations respectively. The actual number of opened stations from  these simulations are shown in Fig.~\ref{fig:rts-v2g2}.  The improvement of load shedding is limited as the optimal solutions of GE-3 to GE-5 in Fig.~\ref{fig:rts-ge} show that to have significant reduction in load shedding we require as many as 25 stations.  In V2G-2 to V2G-4 the increase in unmet battery demands is large when compared to V2G-1. The reason for this is that in V2G-1, the locations of stations were fixed to the optimal solution of the PHEV only system. These locations severely limited the discharge capability of the exchange stations due to the associated line capacities. This result demonstrates that in order for the V2G to benefit both PHEV and the power grid, both models must be considered when determining locations to site exchange stations. In V2G-5, a sufficient number of stations are opened such that both load shedding and unmet battery demands are smaller than in V2G-1. 

In conclusion, at low renewable penetration levels, the V2G system can reduce load shedding caused by renewable generation variability and meet the demand for PHEV battery exchanges. For higher levels of renewables, the benefit of V2G is not as obvious.  In general, the trade-offs between load shedding and unmet battery demands are quite high. These discrepancies may be caused by the relatively low population to load ratio, i.e. there is not enough PHEV demand and only a small number of exchange stations are required.

%Next, we look at a different type of city.

\subsection{Case Study: Miami, Florida}
In our second case, we use Los Alamos National Laboratory's data sets on power grid and transportation networks in Miami, Florida. The transmission power system in the greater Miami area consists of more than 200 buses and over 275 power lines. For this case study, we used peak load data and maximum generation capacity.
%for the buses and the reactance and the line capacity for the lines. 
The overall load is over 6400 MW, and the overall generation and import capacity is greater than 8200 MW.  The second stage problem has 100 scenarios, and in each scenario, loads and generation capacities are independently generated in the same fashion as the first scenario. 
%The load at a bus is a uniform random variable between 0.5 and 1.0 of the maximum load. For a generation bus, if it is a renewable, its generation 
%capacity can be 0, 0.5, or 1 of its maximum capacity equally likely, and otherwise its capacity is same as the maximum capacity. 

In the transportation network, there are over 2500 locations and over 3800 roads. A data set consisting of the gas station locations in the Miami area is used to identify 316 potential locations for exchange stations in the transportation network. 
%and these stations are linked to a nearest bus which is the inserting point for batteries to the power grid. 
%One hundred  pairs of origins and destinations were selected randomly and uniformly from the transportation network. A  shortest path was calculated between each pair of origins and destinations to derive our traffic routes.  
One hundred traffic routes were created by the same method in the last case study.
In 2008, the population of Miami was estimated to be 5,414,712, which yields an estimated 422,348 requests for PHEV batteries based upon the set of ratios described earlier. For each of the 100 scenarios, battery demand of a route is generated randomly and independently with respect to a uniform distribution between 0.5 and 1.5 of the average demand.  Miami represents a city with a high population-to-load ratio when compared to the synthetic city created in the first case study.

\begin{figure}
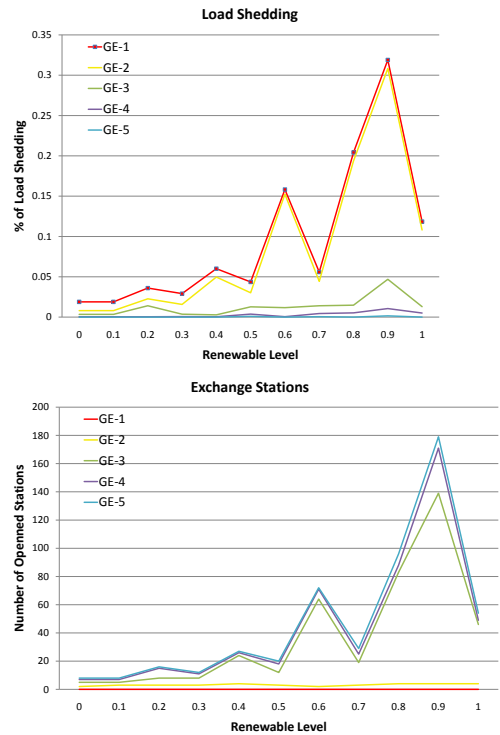

\centering
\psfig{file=PHEV-Figures.pdf,scale=.26,page=6} 
\psfig{file=PHEV-Figures.pdf,scale=.26,page=7}
\caption{\label{fig:miami-ge} \small Miami: load shedding and allocation of exchange stations in generation expansion model with renewables from 0 to 100 percent}
\vspace{-0.6cm}
\end{figure}

Once again, by setting the penalty cost, $h$,  of unmet PHEV demands to zero, we reduce {\bf V2G-2STAGES} to a stochastic generation expansion model. This allows us to understand the maximal benefit V2G could provide to this grid.  We tested this generation expansion model on five different cases GE1-5 in Fig.~\ref{fig:miami-ge}. GE-1 is the base case where no batteries are provided by the exchange stations. In GE-2 to GE-4 we applied different penalty costs on shedding load to investigate the effect of additional exchange stations. With the 11 different renewable penetration levels, the load shedding is not always non-decreasing with respect to the renewable level because of constraints in the power system and the stochastic choice for renewable generation locations. In the results, the location of generators is a key factor in the ability deliver power.  Excess generation capacity at certain locations may not be able to supply loads in the other part of the grid due to capacity constraints.  To elaborate this point, we observed that in GE-1, even without the renewable generation, there was a 2\% load shed even though the overall generation capacity exceeded the overall loads.  As the grid includes more renewables as power sources, the load shedding increases and in the worst case, the load shed reaches 33\%. Increasing the number of exchange station did reduce the load shed. In GE3-GE5, load shedding is reduced to under 5\% for all renewable levels, but the reduction is the result of opening many stations. One interesting case is GE-2 where only 2 to 4 stations are opened and the load shedding is reduced to below 5\% with renewable penetration at 50\% or lower.  For  high  renewable penetration levels, a large number of stations are required to reduce load shed caused by the generation variability. 
%Another reason for  the large reduction in load shed is that  there are more potential locations for exchange stations and it is more likely to find 
%good inserting points for the stations to connect to the grid. 

\begin{figure}
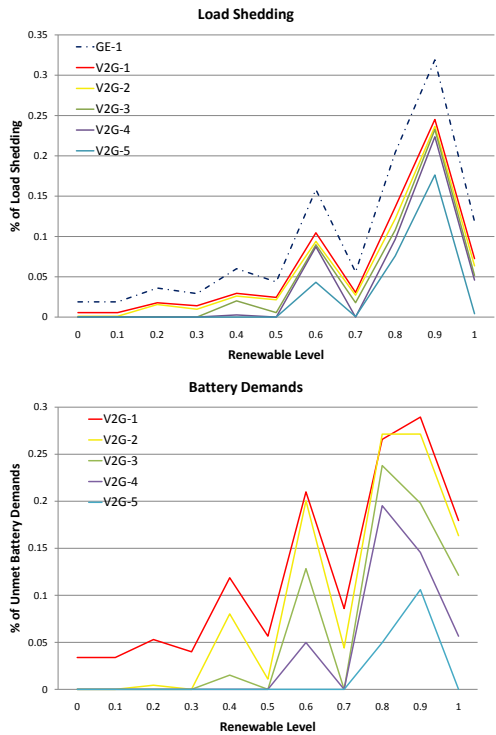

\centering
\psfig{file=PHEV-Figures.pdf,scale=.26,page=8} 
\psfig{file=PHEV-Figures.pdf,scale=.26,page=9}
\caption{\label{fig:miami-v2g1} \small Miami: load shedding and unmet battery demands in {\bf V2G-2STAGES} with renewables from 0 to 100 percent}
\vspace{-0.6cm}
\end{figure}

Using similar steps to the synthetic city case, we removed the power grid from {\bf V2G-2STAGES}, to determine the optimal allocation of the stations (106) to serve all the battery demands. Next, we  evaluated the performance of the full V2G system with respect to  this optimal allocation.  First, we applied the optimal allocation to the power grid and the result is V2G-1 shown in Fig.~\ref{fig:miami-v2g1}. Importantly, the load shed is less than the base case GE-1; however,  the unmet battery demands increased even at low renewable levels. In V2G-2, the total number of opened stations is still set to 106, but {\bf V2G-2STAGES} model can freely choose the locations of stations.  The load shed  is slightly better in V2G-2 than V2G-1 and at the same time the unmet battery demands drop to almost zero for low renewable penetration levels. This indicates that it is important to combine both traffic and grid models when determining the sites to open stations. For V2G-3 to V2G-5, we allowed the number of stations to increase by 10\%, 20\%, and 50\% respectively from 106.  With 10\% increase, V2G-3 shows that {\bf V2G-2STAGES} is able to select locations such that both load shed and unmet battery demands are very low for the renewable penetration level less than and equal to 50\%, indicating that a small amount of additional V2G infrastructure construction can significantly help the power grid.  
%The number of opened stations are shown in Fig.~\ref{fig:miami-v2g2}. It shows that much more stations are needed to improve the performance of the V2G system, and the %improvements are not significant at high levels of renewable penetration. 

In conclusion, with low to medium renewable penetration levels ($\le 50\%$), the V2G system of the greater Miami area  can reduce the variability of renewable sources while maintaining the  service to PHEVs. The location of stations is important for the performance of a V2G system, and integrating both grid and traffic models is important to the planning process. The performance of the V2G system is much better in the Miami case study than in the synthetic city. The difference appears to be a result of  Miami having high population to load ratio, and this high ratio leads to more battery demand, more options for siting exchange stations and greater oppurtunity to benefit the power grid.  
\begin{figure}
\centering
\psfig{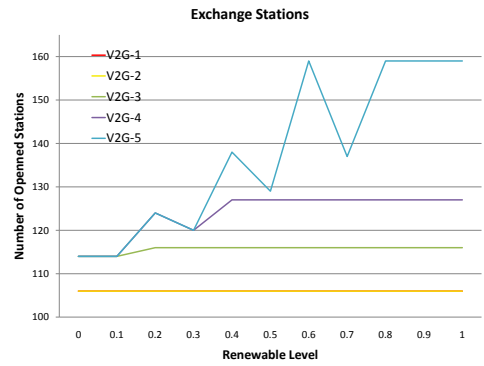} 
\caption{\label{fig:miami-v2g2} \small Miami: allocation of exchange stations in {\bf V2G-2STAGES} with renewables from 0 to 100 percent }
\vspace{-0.6cm}
\end{figure}

\vspace{-.2cm}

\section{Conclusion}\label{sec:conclusion}
In this paper, we investigated the effect of a V2G system on reducing load shed caused by generation variability  of renewables. After introducing a two-stage stochastic programming model to site exchange stations, we applied the model to analyze the V2G system through two case studies.  

\vspace{-.2cm}
 
%\bibliographystyle{plain}

%\bibliography{bib-phev}

\end{document}